\documentclass[12pt]{article}
\usepackage{amssymb}

\begin{document}

\title{On the Semicontinuity in Product Spaces}
\author{Ant\'{o}nio J.B. LOPES-PINTO\thanks{%
Center of Mathematics and Fundamental Applications, University of Lisbon,
Avenida Professor Gama Pinto 2,1649-003, Lisbon Codex, Portugal, E-mail:
lpinto@ptmat.lmc.fc.ul.pt} and Diana Aldea MENDES\thanks{%
Superior Institute of Labour and Management, Avenida das For\c{c}as Armadas,
1699 Lisbon Codex, Portugal}}
\date{}
\maketitle

\begin{abstract}
Let $X,Y$ be topological vector spaces or metric spaces, and let {$f:X\times
Y \rightarrow \Re $} be a real function lower semicontinuous in the first
variable and upper semicontinuous in the second one. It is proved that $f$
is globally measurable. Sierpinski (1925) has been raised this question in
the case $X=Y=\Re $. This particular case was solved by Kempisty (1929). The
actual result has applications in Calculus of Variations.
\end{abstract}

{The classical way to obtain measurability of integrand functions is via the
well-known Carathe\'{o}dory functions. With the development of the Calculus
of Variations and actually, the Shape Optimization ([2,3,4,5]), the weakness
of the above conditions, in particular, the replacement by functions defined
in a topological product with separately semicontinuity is a central
question. Indeed, this is an old and historical question. Sierpinski ([7])
has given an example of a real function defined in the plane which is
separately upper semicontinuous but it isn't measurable. Kempisty ([1]) has
proved that the measurability holds if the function defined in the plane is
lower semicontinuous in one variable and upper semicontinuous in another. In
despite of the motivation of this paper is by some concrete applications,
the authors intend this as a technical note in order to clarify certain
situations for real functions defined in topological spaces.}

{In this note, using the method of Kempisty ([1]), we obtain a positive
result for real functions defined in product of metric spaces and we point
out some difficulties for the extension to arbitrary product topologies. The
idea base of the proof consists at the construction of upper semicontinuous
functions and lower semicontinuous functions which approach arbitrarily the
function from below and from above, respectively, in order to obtain a lower
semicontinuous multivalued map. Therefore, this kind of proof consists of
the following steps: 1) to construct \ adequate lower and upper
semicontinuous functions strictly related with the given function; 2) to
define functions of this type which are enough close to the given function;
3) to use some continuous selection theorem for multis with real closed
interval values. The last one suggests an application of classical selection
theorems for lower semicontinuous multivalued maps.}

{The paper is organized in order to put in evidence the aspects cited above.
We star showing that 1) holds for topological vector spaces but 2) doesn't
hold and it isn't available an adequate continuous selection theorem.
Summing up the results, it is expectable that an alternative proof works out
to uniform spaces (at least, for paracompact spaces, once every open cover
is even). Actually, we prove essentially a little more: every real function
defined on a metric product space which is lower semicontinuous in the first
variable and upper semicontinuous in the second one \ is a second category
function of Baire at most, and hence, a Borel measurable function.}

\section{Lower and upper approximations}

{Let $X,Y$ be two nonempty topological spaces and let $f:X\times
Y\rightarrow \,\Re \,$be a real valued function.}

{Let $V$ be a neighborhood of $y_{0}$ in $Y$. We denote by $M_{V}^{2}\left(
x_{0},y_{0}\right) $ the upper bound of the function $f\left( x_{0},y\right)
\,$on the neighborhood $V$ and by $m_{V}^{2}\left( x_{0},y_{0}\right) \,$the
lower bound of the function $f\left( x_{0},y\right) \,$on $V$, i.e., 
\[
M_{V}^{2}\left( x_{0},y_{0}\right) =\sup_{y\in V}f\left( x_{0},y\right) 
\]
} 
\[
m_{V}^{2}\left( x_{0},y_{0}\right) =\inf_{y\in V}f\left( x_{0},y\right) .
\]
\textbf{Lemma 1.1}{\ }\textit{Let }$X$\textit{\ be a nonempty topological
space}$,Y\,$\textit{be a topological vectorial space and }$f:X\times
Y\rightarrow \,\Re $\textit{\ be a real valued function. Let fix an open
neighborhood }$\mathcal{W}_{0}\,$\textit{of }$0\in Y$\textit{. If }$f\left(
\cdot ,y\right) $\textit{\ is lower semicontinuous at }$x_{0}\,\,$\textit{%
for any }$y\in Y$\textit{\ then the upper bound } 
\[
M_{\mathcal{W}_{0}}^{2}\left( x,y\right) :=\sup_{z\in y+\mathcal{W}%
_{0}}f\left( x,z\right) 
\]
\textit{is lower semicontinuous at }$\left( x_{0},y_{0}\right) $\textit{\ \
for any }$y_{0}\in Y$\textit{.}

\ 

\textbf{Proof}{\ \ Let $\left( x_{0},y_{0}\right) \in X\times Y$.}

{We are interested to prove that \ for every $L\in \Re $ such that $L<M_{%
\mathcal{W}_{0}}^{2}\left( x_{0},y_{0}\right) $ (eventually, $M_{\mathcal{W}%
_{0}}^{2}\left( x_{0},y_{0}\right) =+\infty $), there exist a neighborhood $%
\mathcal{V}$ of $x_{0}$ and a neighborhood $\mathcal{W}$ of $0$ such that
for every $\left( x,y\right) \in \mathcal{V}\times \left( y_{0}+\mathcal{W}%
\right) $ it is $L<$ $M_{\mathcal{W}_{0}}^{2}\left( x,y\right) $.}

{Let $y_{1}\in y_{0}+\mathcal{W}_{0}$ be such that 
\[
L<f\left( x_{0},y_{1}\right) . 
\]
As $f\left( \cdot ,y_{1}\right) $ is lower semicontinuous at $x_{0}$, there
exists a neighborhood $\mathcal{V}$ of $x_{0}$ such that 
\begin{equation}
L<f\left( x,y_{1}\right) \hspace*{0.2in}\mathnormal{\mathnormal{for every }}%
x\in \mathcal{V}.\mathnormal{\ }  \label{1}
\end{equation}
If we choose a balanced neighborhood $\mathcal{W}\,$ of \ $0\in Y$ \ such
that \ 
\[
y_{1}-y_{0}+\mathcal{W}\subseteq \mathcal{W}_{0}, 
\]
it follows 
\[
L<\mathnormal{\ }M_{\mathcal{W}_{0}}^{2}\left( x,y\right) \hspace*{0.2in}%
\mathnormal{\mathnormal{\ for every }}\left( x,y\right) \in \mathcal{V}%
\times \left( y_{0}+\mathcal{W}\right) . 
\]
Indeed, if $y\in y_{0}+\mathcal{W},$ it is 
\[
y_{1}-y=\left( y_{1}-y_{0}\right) +\left( y_{0}-y\right) \in \left(
y_{1}-y_{0}\right) +\mathcal{W}\subseteq \mathcal{W}_{0}\mathnormal{\ } 
\]
Hence, $y_{1}\in y+\mathcal{W}_{0}$ and according to $\left( \ref{1}\right) $
\ 
\[
L<f\left( x,y_{1}\right) \leq M_{\mathcal{W}_{0}}^{2}\left( x,y\right)
\hspace*{0.2in}\mathnormal{\mathnormal{\ for every }}\left( x,y\right) \in 
\mathcal{V}\times \left( y_{0}+\mathcal{W}\right) \mathnormal{\ .} 
\]
}

{\hfill }$\square $

\textbf{Corollary 1.2 }\textit{Let }$X$\textit{\ be a nonempty topological
space}$,Y\,$\textit{be a topological vectorial space and let }$f:X\times
Y\rightarrow \,\Re $\textit{\ be a real valued function. We fix an open}$\,$%
\textit{neighborhood }$V_{0}\,$\textit{of }$0$\textit{\ }$\left( \mathnormal{%
zero}\right) $\textit{. If }$f\left( \cdot ,y\right) $\textit{\ is upper
semicontinuous at }$x_{0}\,$\textit{for any }$y\in Y$\textit{\ then the
lower bound } 
\[
m_{V_{0}}^{2}\left( x,y\right) =\inf_{z\in y+V_{0}}f\left( x,z\right) 
\]
\textit{\ is upper semicontinuous at }$\left( x_{0},y_{0}\right) $\textit{\
for any }$y_{0}\in Y$\textit{.}

\textbf{Proof} It is enough to replace in the lemma $f$ by $-f$. $\square $

\textbf{Remark}{\ We consider now 
\[
M_{V_{0}}^{2}\left( x_{0},y_{0}\right) \downarrow M^{2}\left(
x_{0},y_{0}\right) \geq f\left( x_{0},y_{0}\right) 
\]
} i.e., $M^{2}\,$it is the pointwise limit of a nonincreasing sequence of
lower semicontinuous functions, or $M^{2}\,$is a function of type $ul$ via
M. Young.

{Analogous 
\[
m_{V_{0}}^{2}\left( x_{0},y_{0}\right) \uparrow m^{2}\left(
x_{0},y_{0}\right) \leq f\left( x_{0},y_{0}\right) 
\]
} i.e., $m^{2}$ is the pointwise limit of a nondecreasing sequence of upper
semicontinuous functions, or$\,m^{2}$ is a function of type $lu$ via M.
Young.

{An analogous lemma holds for a metric space.}

{Let $\alpha >0$ and \textbf{B}$\left( y_{_{0}};\alpha \right) :=\left\{
y\in Y:d\left( y,y_{0}\right) <\alpha \right\} $ be the ball with center at $%
y_{0}\in Y$ and radius $\alpha $ in the metric space $Y$. We denote by $%
M_{\alpha }^{2}\left( x_{0},y_{0}\right) $ the upper bound of the function $%
f\left( x_{0},y\right) \,$on the ball \textbf{B}$\left( y_{_{0}};\alpha
\right) $ and by $m_{\alpha }^{2}\left( x_{0},y_{0}\right) \,$the lower
bound of the function $f\left( x_{0},y\right) \,$on the ball \textbf{B}$%
\left( y_{_{0}};\alpha \right) $, i.e., 
\[
M_{\alpha }^{2}\left( x_{0},y_{0}\right) =\sup_{y\in \mathnormal{B}\left(
y_{_{0}};\alpha \right) }f\left( x_{0},y\right) 
\]
} 
\[
m_{\alpha }^{2}\left( x_{0},y_{0}\right) =\inf_{y\in \mathnormal{B}\left(
y_{_{0}};\alpha \right) }f\left( x_{0},y\right) . 
\]

\textbf{Lemma 1.3}{\ }\textit{Let }$X$\textit{\ be a nonempty topological
space}$,\left( Y,d\right) \,$\textit{be a nonempty metric space and let }$%
f:X\times Y\rightarrow \,\Re $\textit{\ be a real valued function. Let fix }$%
\alpha >0$\textit{. If }$f\left( \cdot ,y\right) $\textit{\ is lower
semicontinuous at }$x_{0}\,$\textit{\ for any }$y\in Y$\textit{\ then the
upper bound } 
\[
M_{\alpha }^{2}\left( x,y\right) :=\sup_{z\in \mathnormal{B}\left( y;\alpha
\right) }f\left( x,z\right) 
\]
\textit{is lower semicontinuous at }$\left( x_{0},y_{0}\right) $\textit{\
for any }$y_{0}\in Y$\textit{.}

\textbf{Proof}{. Let $\left( x_{0},y_{0}\right) \in X\times Y$. }

{We are interested to prove that \ for every $L\in \Re $ such that $%
L<M_{\alpha }^{2}\left( x_{0},y_{0}\right) $ (eventually, $M_{\alpha
}^{2}\left( x_{0},y_{0}\right) =+\infty $), there exist a neighborhood $%
\mathcal{V}$ of $x_{0}$ and $r>0$ such that for every $\left( x,y\right) \in 
\mathcal{V}\times $\textbf{B}$\left( y_{0};r\right) $ it is $L<$ $M_{\alpha
}^{2}\left( x,y\right) $.}

{Let $y_{1}\in $\textbf{\ B}$\left( y_{0};\alpha \right) $ be such that 
\[
L<f\left( x_{0},y_{1}\right) .
\]
As $f\left( \cdot ,y_{1}\right) $ is lower semicontinuous at $x_{0}$, there
exists a neighborhood $\mathcal{V}$ of $x_{0}$ such that 
\begin{equation}
L<f\left( x,y_{1}\right) \hspace*{0.2in}\mathnormal{forevery}x\in \mathcal{V}%
.\mathnormal{\ }  \label{2}
\end{equation}
If we choose $r:0<r<\alpha -d\left( y_{0},y_{1}\right) $ it follows 
\[
L<\mathnormal{\ }M_{\alpha }^{2}\left( x,y\right) \hspace*{0.2in}\mathnormal{%
forevery}\left( x,y\right) \in \mathcal{V}\times \mathnormal{\mathbf{\ B}}%
\left( y_{0};r\right) .
\]
Indeed, if $y\in $B$\left( y_{0};r\right) $ it is $d\left( y,y_{1}\right)
\leq d\left( y,y_{0}\right) +d\left( y_{0},y_{1}\right) <\alpha $. Hence, $%
y_{1}\in $B$\left( y;\alpha \right) $ and, according to (\ref{2})
\[
L<f\left( x,y_{1}\right) \leq M_{\alpha }^{2}\left( x,y\right)
\hspace*{0.2in}\mathnormal{forevery}\left( x,y\right) \in \mathcal{V}\times 
\mathnormal{B}\left( y_{0};r\right) \mathnormal{\ .}
\]
}

\hfill ${\square }$

\textbf{Corollary 1.4 }\textit{Let }$X$\textit{\ be a nonempty topological
space, }$\left( Y,d\right) $\textit{\ be a nonempty metric space and let }$%
f:X\times Y\rightarrow \,\Re $\textit{\ be a real valued function. We fix }$%
\alpha >0$\textit{. If }$f\left( \cdot ,y\right) $\textit{\ is upper
semicontinuous at }$x_{0}$\textit{\ for any y}$\in Y$\textit{\ then the
lower bound}

\[
m_{\alpha }^{2}\left( x,y\right) :=\inf_{z\in B\left( y,\alpha \right)
}f\left( x,z\right) 
\]
\textit{is upper semicontinuous at }$\left( x_{0},y_{0}\right) ,\,$\textit{%
for any }$y_{0}\in Y$\textit{.}

\textbf{Proof}{. It is enough to replaced in the preceding lemma $f$ by $-f.$
$\square $}

\textbf{Corollary 1.5 }{\ }\textit{Let }$X$\textit{\ be a nonempty
topological space}$,Y\,$\textit{be a Fr\'{e}chet (or a normed) space and let 
}$f:X\times Y\rightarrow \,\Re $\textit{\ be a real valued function. We fix }%
$\alpha >0$\textit{. If }$f\left( \cdot ,y\right) $\textit{\ is lower
semicontinuous at }$x_{0}\,$\textit{for any }$y\in Y$\textit{\ then the
upper bound } 
\[
M_{\alpha }^{2}\left( x,y\right) :=\sup_{z\in \mathnormal{B}\left( y;\alpha
\right) }f\left( x,z\right) 
\]
\textit{is lower semicontinuous at }$\left( x_{0},y_{0}\right) $\textit{,
for any }$y_{0}\in Y$\textit{. }

\section{Measurability of the function.}

We recall that if $X$ is a topological space and \ \textnormal{$%
a,b,X\rightarrow $ $\Re $} are upper semicontinuous and lower semicontinuous
at $x_{0},$ respectively, and $a\leq b,$ then the multivalued map 
\[
\Gamma :X\rightarrow 2^{\mathnormal{\Re }},x\mapsto [a\left( x\right)
,b\left( x\right) ]
\]
is lower semicontinuous at $x_{0}.\,$The Michael's Selection theorem ([6])
asserts that if $X$ is paracomapct (in particular, if $X$ is a metric space)
then there exists a continuous selection (actually, this theorem is an
extension of a Hahn theorem ([1]))

\textbf{Theorem 2.1 }{\ \textit{Let }$X,Y$\textit{\ be two nonempty metric
spaces and let }$f:X\times Y\rightarrow \,\Re $\textit{\ be a bounded real
valued function.}$\,$\textit{If }$f\left( x,y\right) \,$\textit{it is lower
semicontinuous in relation to one of the variables and it is upper
semicontinuous in relation to another one, then }$f\left( x,y\right) \,$%
\textit{is of the first category of Baire at most}.}

\textbf{Proof}{. We assume that $f\left( \cdot ,y\right) \,$is lower
semicontinuous for every $y$ and that $f\left( x,\cdot \right) $ \ is upper
semicontinuous for every $x$. }

{Let define for every $n{\in \mathbb{N}}$ 
\[
m_{\frac{1}{n}}^{1}\left( x,y\right) :=\inf_{w\in B\left( x;\frac{1}{n}%
\right) }f\left( w,y\right) ,\hspace*{0.2in}M_{\frac{1}{n}}^{2}\left(
x,y\right) :=\sup_{z\in \mathnormal{B}\left( y;\frac{1}{n}\right) }f\left(
x,z\right) 
\]
}

{It is clear that for every $\left( x,y\right) \in X\times Y$ the sequences $%
\left( m_{\frac{1}{n}}^{1}\left( x,y\right) \right) _{n}$ and $\left( M_{%
\frac{1}{n}}^{2}\left( x,y\right) \right) _{n}$ are nondecreasing and
nonincreasing, respectively, say, 
\[
m_{\frac{1}{n}}^{1}\left( x,y\right) {\uparrow }_{n} m\left(
x,y\right) \mathnormal{\ \hspace*{0.2in}and\hspace*{0.2in}}M_{\frac{1}{n}%
}^{2}\left( x,y\right)  {\downarrow }_{n}  M\left( x,y\right) \,. 
\]
By other side,}

\[
m_{\frac{1}{n}}^{1}\left( x,y\right) \leq f\left( x,y\right) \leq M_{\frac{1%
}{n}}^{2}\left( x,y\right) ,\hspace*{0.2in}\left( x,y\right) \in X\times
Y,\hspace*{0.1in}n{\in \mathbb{N}}\mathnormal{\ .} 
\]
We claim that \ 
\[
m\left( x,y\right) =f\left( x,y\right) =M\left( x,y\right)
,\hspace*{0.2in}\left( x,y\right) \in X\times Y\mathnormal{\ .} 
\]
Indeed, let assume that for some $\left( x_{0},y_{0}\right) $ 
\[
m\left( x_{0},y_{0}\right) <f\left( x_{0},y_{0}\right) . 
\]
Let $\epsilon :0<\epsilon <f\left( x_{0},y_{0}\right) -m\left(
x_{0},y_{0}\right) .\,$

As $f\left( \cdot ,y_{0}\right) $ is lower semicontinuous at $x_{0},$ it
exists $n_{0}{\in \mathbb{N}}$ such that 
\[
f\left( x_{0},y_{0}\right) <f\left( w,y_{0}\right) +\frac{\epsilon }{2}%
,\hspace*{0.2in}w\in B\left( x_{0};\frac{1}{n_{0}}\right) . 
\]
Hence, 
\[
m_{\frac{1}{n_{0}}}^{1}\left( x_{0},y_{0}\right) \leq m\left(
x_{0},y_{0}\right) <f\left( x_{0},y_{0}\right) <f\left( w,y_{0}\right)
+\epsilon ,\hspace*{0.2in}w\in B\left( x_{0};\frac{1}{n_{0}}\right) . 
\]
But, by definition of $m_{\frac{1}{n_{0}}}^{1}\left( x_{0},y_{0}\right) $,
it exists $w_{0}\in B\left( x_{0};\frac{1}{n_{0}}\right) $ such that $%
f\left( w_{0},y_{0}\right) +\epsilon <f\left( x_{0},y_{0}\right) $ \ and the
contradiction \ 
\[
f\left( w_{0},y_{0}\right) +\epsilon <f\left( x_{0},y_{0}\right) <f\left(
w_{0},y_{0}\right) +\frac{\epsilon }{2} 
\]
\ \ Similarly, it is $f\left( x_{0},y_{0}\right) =M\left( x_{0},y_{0}\right) 
$.

{Now we construt a sequence of continuous functions which converges
pointwise to $f.\,$}

Indeed, from Hahn Theorem, for every $n{\in \mathbb{N}}$ there exists an
intermediary continuous function $f_{n}$ between $m_{\frac{1}{n}}^{1}$and $%
M_{\frac{1}{n}}^{2}$, \ i.e., \ \ 

\[
m_{\frac{1}{n}}^1\left(x,y\right)\leq f_n\left(x,y\right)\leq M_{\frac{1}{n}%
}^2\left(x,y\right),\hspace*{0.2in}\left(x,y\right)\in X\times Y, 
\]
as $M_{\frac{1}{n}}^2$ and $m_{\frac{1}{n}}^1$ are lower semicontinuous and
upper semicontinuous, respectively (lemma 1.3. and corollary 1.4.).

{So, $f\left( x,y\right) $ is a function of the first class of Baire (can be
represented as the limit of an everywhere convergent sequence of continuous
functions).}$\square $

\textbf{Corollary 2.2}{\ \textit{Let }$X,Y$\textit{\ be two nonempty metric
spaces and let }$f:X\times Y\rightarrow \,\Re $\textit{\ be a real valued
function.}$\,$\textit{If }$f\left( x,y\right) $\textit{\ }$\,$\textit{is
lower semicontinuous in relation to one of the variables and it is upper
semicontinuous in relation to another one, then }$f\left( x,y\right) \,$%
\textit{is of the second category of Baire at most}.}

\textbf{Proof}{. Let $\alpha >0$ and let define $f_{\alpha }:X\times
Y\rightarrow \,\Re $, 
\[
\left( x,y\right) \mapsto \left\{ 
\begin{array}{lll}
f\left( x,y\right) & \mathnormal{if } & \left| f\left( x,y\right) \right|
\leq \alpha \\ 
\alpha & \mathnormal{if} & f\left( x,y\right) >\alpha \\ 
-\alpha & \mathnormal{if} & f\left( x,y\right) <-\alpha
\end{array}
\right\} . 
\]
$f_{\alpha }$ has the same type of semicontinuity as $f$. So $f$ is the
pointwise limit of the sequence of bounded functions $f_{\frac{1}{n}}$ with
the same type of semicontinuity as $f$. According to the preceding theorem,
every $f_{\frac{1}{n}}$ is $\,$is of the first category of Baire. Hence, $f$
is a function of the second category of Baire. }$\square $

\textbf{Theorem 2.3}{\ }\textit{Let }$X,Y$\textit{\ be two nonempty metric
spaces and let }$f:X\times Y\rightarrow \,\Re $\textit{\ be a real valued
function.}$\,$\textit{If }$f\left( x,y\right) \,$\textit{it is lower
semicontinuous in relation to one of the variables and it is upper
semicontinuous in relation to another one, then }$f\left( x,y\right) \,$%
\textit{is Borel measurable.}

\textbf{Proof}{\ The result follows immediately as $f$ is a function of the
second category of Baire, and every Baire function is Borel measurable. }$%
\square ${\ }

{\LARGE References}

$\left[ 1\right] $ S. Kempisty, Sur les fonctions semicontinues par rapport
\`{a} chacune de deux variables, Fund. Math. XIV (1929), 237-241.

$\left[2\right]$ A. Lopes-Pinto \& Diana Mendes, On the Semicontinuity and
the Measurability in ODE and Variational Principles. In preparation. \ 

$\left[ 3\right] $ P. Marcellini, Alcune osservazioni sull'esistenza del
minimo di integrali del calcolo delle variazioni senza ipotesi di
convessit\'{a}, Rendiconti Mat., 13 (1980), 271-281. \ 

$\left[4\right]$ P. Marcellini \& C. Sbordone, Semicontinuity Problems in
the Calculus of Variations, Nonlinear Analysis, Theory, Methods \&
Applications, 4,2 (1980), 241-257.

$\left[5\right]$ P. Marcellini, A relation between existence of minima for
non convex integrals and uniqueness for non strictly convex integrals of the
calculus of variations, Lecture Notes in Mathematics, Vol. 979, pp. 216-232.
Springer Verlag, New York (1983).

$\left[6\right]$ E. Michael, Continuous selections I, Annals of Mathematics,
63, 2 (1956), 361-382.

$\left[7\right]$ W. Sierpinski, Funkcje przedstawialne analityczne, Fund.
Math. \ (1925), 68.

\end{document}